\documentclass[11pt]{amsart}

\usepackage{amsmath,amssymb,epsf, epic, eepic}
\usepackage{psfig}
\usepackage{amsfonts}
\usepackage{times}

\input xy
\xyoption{all}

\newtheorem{thm}{Theorem}[section]

\newtheorem{cor}[thm]{Corollary}
\newtheorem{prop}[thm]{Proposition}

\theoremstyle{definition}

\newcommand{\bQ}{\mathbb{Q}}

\newcommand{\ra}{\rightarrow}

\newcommand{\spaces}{\;\;\;\;\;\;\;}
\newcommand{\ot}{\otimes}
\newcommand{\ol}{\widetilde}
\newcommand{\ob}{\overline}
\newcommand{\cc}[3]{C^{{#1},{#2}}(#3)}
\newcommand{\ccol}[3]{\ob{C}^{{#1},{#2}}(#3)}
\newcommand{\kh}[3]{K\!H^{{#1},{#2}}(#3)}
\newcommand{\khr}[4]{K\!H^{{#1},{#2}}(#3;#4)}
\newcommand{\lee}[2]{Lee^{{#1}}(#2)}
\newcommand{\dk}[1]{D_{({#1})}}
\newcommand{\dkb}[1]{\widetilde{D}_{({#1})}}

\newcommand{\lk}{lk}

\title{A spectral sequence for Khovanov homology with an application to 
$(3,q)$-torus links}
\author{Paul Turner}
\address{School of Mathematical and Computer Sciences \\Heriot-Watt
  University\\ Edinburgh EH14 4AS\\Scotland}
\email{paul@ma.hw.ac.uk}

\begin{document}


\begin{abstract}
A spectral sequence converging to Khovanov homology is constructed which is applied to calculate the rational Khovanov homology of $(3,q)$-torus links.
\end{abstract}

\maketitle

\section{Introduction}
There is a lack of theoretical computational
tools for Khovanov homology when compared to, say, the homology of spaces where one has a range
of long exact sequences and spectral sequences at hand. There
is one long exact sequence used in Khovanov homology, 
the {\em skein exact sequence}, implicit in \cite{khovanov} and 
explicit in \cite{viro} which is formed as follows. A given 
crossing of an  link diagram $D$ can be resolved in two ways: to the 0-smoothing  giving a new
diagram $D^\prime$ and to the 1-smoothing giving  a new diagram
$D^{\prime\prime}$. There is then a short exact sequence of Khovanov complexes
\[
0 \ra C(D^{\prime\prime}) \ra C(D)\ra C(D^\prime) \ra 0
\]
which gives rise to a long exact sequence in homology. One can
repeatedly apply this long exact sequence, but it requires careful
book keeping. This is essentially what leads to the spectral
sequence defined in this paper.

We start with a collection of $m$ crossings of a diagram $D$. For
$1\leq k \leq m$ let $\dk k$ be the diagram obtained from $D$ by
resolving the crossings $1,\ldots, k$ to 1-smoothings and let $\dkb k$
be the diagram obtained from $D$ by resolving the crossings $1,\ldots,
k-1$ to 1-smoothings and crossing $k$ to a 0-smoothing. The idea is
that the diagrams $\dkb k$ and $\dk m$ might be simpler to
handle than the diagram $D$. By defining appropriate constants
$a_k,b_k,\ol{a}_k,\ol{b}_k, A_k,$ and $ B_k$ we arrive at the following
result where $j$ is fixed.

\vspace*{4mm}

{\bf Proposition \ref{prop:ss}} There is a spectral sequence $(E_r^{*,*}, d_r\colon E_r^{s,t} \ra E_r^{s+r,t-r+1})$ converging to 
$\kh * j D$ with $E_1$-page given by
\[
E_1^{s,t} = \begin{cases}
\kh {s+t+A_s+\ol{a}_{s+1}} {j+B_s+\ol{b}_{s+1}} {\dkb {s+1}} & s=0,\ldots,m-1\\
\kh {m+t+A_m} {j+B_m} {\dk m} & s=m\\
0 & s<0 \text{ or } s> m
	    \end{cases}
\]
 
\vspace*{4mm}

As an application of this spectral sequence we compute the
rational Khovanov homology of $(3,q)$-torus links. It is easy to guess
what the result is, based on available computer calculations, but by
combining the above spectral sequence with another spectral sequence
(Lee's spectral sequence) we prove the result for all $q$.

\section{The spectral sequence}
Let $R$ be a commutative ring with unit and let $D$ be an oriented
link diagram with $n$ crossings. As is now familiar one can construct
the Khovanov complex by placing the $2^n$ smoothings of $D$ on the
vertices of the cube $\{0,1\}^n$. To each smoothing $\alpha$ one then
assigns the $R$-module $V_\alpha = V^{\ot k_\alpha}\{r_\alpha\}$ where
$k_\alpha$ is the number of circles in the smoothing, $r_\alpha$ is
the number of 1's in $\alpha$ and shifts are defined by $(W\{l\})^m =
W^{m-l}$. The module $V$ is the graded, rank two, free $R$-module with
generators $1$ and $x$ in degree $1$ and $-1$ respectively. The
underlying module of the {\em unnormalised} Khovanov complex
$\ccol ** D$ is defined by
\[
\ccol i * D = \bigoplus_{\substack{\alpha\in \{0,1\}^n\\ r_\alpha =i}}
V_\alpha.
\]
The construction of the differential is by now well known (see, for
example, \cite{khovanov} or  \cite{barnatan}) and uses a Frobenius
algebra structure on $V$.

Bi-graded complexes may be shifted in each of the degrees and for a
bi-graded module $W^{**}$ we define
\[
(W^{*,*}[l]\{m\})^{i,j} = W^{i-l,j-m}.
\]

Suppose $D$ has $n^+$ positive crossings and $n^-$ negative crossings, then the {\em normalised} Khovanov complex $\cc ** D$ is defined by 
\[
\cc i j D = (\ccol ** D [-n^-]\{n^+-2n^-\})^{i,j}
\]
and the Khovanov homology of $D$ is defined as the homology of this complex.

Now let us consider a collection of $m$ crossings of the diagram $D$
and number these $1,\ldots, m$. For $k=1,\ldots,m$ let $\dk k$ be the
diagram obtained from $D$ by resolving the crossings $1,\ldots, k$ to
1-smoothings and let $\dkb k$ be the diagram obtained from $D$ by
resolving the crossings $1,\ldots, k-1$ to 1-smoothings and crossing
$k$ to a 0-smoothing. We also define $\dk 0$ and $\dkb 0$ to be the
original diagram $D$.

There is a decomposition of modules
\[
\ccol i j {\dk {k-1}} = \ccol i j {\dkb {k}} \oplus \ccol {i-1}{j-1} {\dk {k}}
\]
and in fact $\ccol {*-1}{*-1} {\dk {k}}$ is a sub-complex of $\ccol ** {\dk {k-1}}$. Thus, there is a short exact sequence
\begin{equation}\label{eq:ses}
\xymatrix{
0 \ar[r] & \ccol {*}{*} {\dk {k}}[1]\{1\} \ar[r] & \ccol ** {\dk {k-1}} \ar[r] &
\ccol ** {\dkb {k}} \ar[r] & 0.
}
\end{equation}

This is just the usual short exact sequence giving the long exact
sequence mentioned in the introduction for the diagram $\dk k$
resolving the $k$'th crossing in our set
of $m$ crossings.

We now discuss orientations for the diagrams $\dk k$ and $\dkb
k$. Suppose that we already have an orientation for $\dk {k-1}$. If the $k$'th
crossing is positive then $\dkb {k}$ inherits an orientation because
for positive crossings the 0-smoothing is the oriented
resolution. There is no orientation of $\dk k$ consistent with the
orientation of $\dk {k-1}$ so choose any orientation for $\dk
k$. Similarly if the $k$'th crossing is negative then $\dk {k}$
inherits an orientation and we choose any orientation for $\dkb k$.
The diagram $D=\dk 0= \dkb 0$ comes with an orientation so the process above has
somewhere to start.

Now for $k=0,\ldots , m$ define
\begin{eqnarray*}
n_k^+  & = & \text{number of positive crossings in }\dk k\\
n_k^-  & = & \text{number of negative crossings in }\dk k\\
\ol{n}_k^+  & = & \text{number of positive crossings in }\dkb k\\
\ol{n}_k^-  & = & \text{number of negative crossings in }\dkb k\\
\end{eqnarray*}
We define additional constants associated to $\dk k$
and $\dkb k$ as follows.

If the $k$'th crossing is positive then set

\[
a_k  =  n_{k-1}^- - n_k^- -1 \;\;\;\text{ and }\;\;\;
\ol{a}_k  =  0.
\]

If the $k$'th crossing is negative then set
\[
a_k  =  0\;\;\;\text{ and }\;\;\;
\ol{a}_k  =   n_{k-1}^-  -  \ol{n}_k^-.
\]
For convenience we also define (for positive and negative crossings)
\[
b_k = 3a_k +1 \;\;\;\text{ and }\;\;\; \ol{b}_k = 3\ol{a}_k -1.
\]

These constants help us to
write down the short exact sequence (\ref{eq:ses}) in terms of
normalised Khovanov homology. 

\begin{prop}\label{prop:ses}
For each $k=1,\ldots, m$ there is a short exact sequence of complexes
{\small
\begin{equation*}
\xymatrix{
0 \ar[r] & \cc {*}{*} {\dk {k}}[-a_k]\{-b_k\} \ar[r] & \cc ** {\dk {k-1}} \ar[r] &
\cc ** {\dkb {k}}[-\ol{a}_k]\{-\ol{b}_k\} \ar[r] & 0
}
\end{equation*}
}
\end{prop}

\begin{proof}
We shift the entire sequence  (\ref{eq:ses}) by $[-n_{k-1}^-]\{n_{k-1}^+ - 2n_{k-1}^-\}$. One can then readily verify (treating positive and negative crossings separately) that  
\begin{eqnarray*}
-n_{k-1}^- +1 & = & -n_k^- - a_k, \\
-n_{k-1}^- & = & -\ol{n}_k^- -\ol{a}_k,\\
n_{k-1}^+ - 2n_{k-1}^- + 1 & = & n_{k}^+ - 2n_{k}^- -b_k, \\
n_{k-1}^+ - 2n_{k-1}^- &=&   \ol{n}_{k}^+ - 2\ol{n}_{k}^- - \ol{b}_k.  
\end{eqnarray*}
\end{proof}

We now define
\[
A_k  =  \sum_{i=1}^k a_i, \;\; \text{and} \;\; B_k =  \sum_{i=1}^k b_i = 3A_k + k.
\]
and set $A_0=B_0=0$.

From now on we fix $j$. We define a filtration on $\cc * j D$ by
\[
F^k\cc * j D = \cc *j {\dk k}[-A_k]\{-B_k\} \;\;\;\;\; k=0,\ldots,m.
\]

It follows immediately from Proposition \ref{prop:ses} that $F^k\cc *
j D \subset F^{k-1}\cc * j D$ and for $k>m$ we set $F^k\cc * j D=0$. It
is clear that the filtration is bounded and so there is an associated
spectral sequence.

\begin{prop}\label{prop:ss}
There is a spectral sequence  $(E_r^{*,*}, d_r\colon E_r^{s,t} \ra E_r^{s+r,t-r+1})$ converging to 
$\khr * j D R$ with $E_1$-page given by
\[
E_1^{s,t} = \begin{cases}
\khr {s+t+A_s+\ol{a}_{s+1}} {j+B_s+\ol{b}_{s+1}} {\dkb {s+1}} R & s=0,\ldots,m-1\\
\khr {m+t+A_m} {j+B_m} {\dk m} R & s=m\\
0 & s<0 \text{ or } s> m
	    \end{cases}
\]
\end{prop}

\begin{proof}
By using the filtration above there is a spectral sequence with
\[
E_0^{s,t} = \frac{ F^{s}\cc {s+t} j D}{ F^{s+1}\cc {s+t} j D}.
\]
By applying Proposition \ref{prop:ses} there is a short exact sequence
for $0\leq s < m$ as follows.  
{\small
\[
\xymatrix{
0 \ar[r] &  F^{s+1}\cc {s+t} j D \ar[r] & F^{s}\cc {s+t} j D  \ar[r] &
\cc {s+t+A_s+\ol{a}_{s+1}} {j+B_s+\ol{b}_{s+1}}  {\dkb {s+1}} \ar[r] & 0
}
\]
}
Since this is a short exact sequence of complexes the differential
$d_0$, which is induced by the differential on $F^s\cc{*}j D$, can be
identified with the differential on the right hand side, that is, in  the
complex $\cc ** {\dkb {s+1}}$. In particular the homology of
$E_0^{*,*}$ is given by the homology (in suitable gradings) of $\cc ** {\dkb {s+1}}$, namely the Khovanov homology of $\dkb
{s+1}$.

When $s=m$ we have $E_0^{m,t} = F^m \cc {m+t} j D = \cc {m+t +A_m}
{j+B_m} {\dk m}$ and so $E_1$ is again as claimed.
\end{proof}

The differential $d_1$ on the $E_1$-page can be understood as follows.
There is a decomposition (of modules)
\begin{align*}
\cc ** {\dk s} = & \cc {*+\ol{a}_{s+1}}{* + \ol{b}_{s+1}}  {\dkb {s+1}} \\
& \spaces\spaces \oplus \cc {*+a_{s+1}  +  \ol{a}_{s+2}}{*+b_{s+1} + \ol{b}_{s+2}}{\dkb {s+2}}\\
& \spaces\spaces\spaces\spaces \oplus  \cc {*+a_{s+1} + a_{s+2}}{*+b_{s+1}+ b_{s+2}}{\dk {s+2}}
\end{align*}
and with respect to this the  differential on $\cc ** {\dk s}$ can be written as a matrix 
\[
\begin{pmatrix}\ol{\delta}_{s+1} & 0 & 0\\ 
\delta & \ol{\delta}_{s+2}& 0\\
\delta^\prime & \delta^{\prime\prime} & {\delta}_{s+2} \end{pmatrix}.
\]

The differential on the $E_1$-page of the spectral sequence is the map 
\[
\delta\colon \cc {*+\ol{a}_{s+1}}{* + \ol{b}_{s+1}}{\dkb {s+1}} \ra \cc {*+a_{s+1}  +  \ol{a}_{s+2}}{*+b_{s+1} + \ol{b}_{s+2}}{\dkb {s+2}}
\]
in the above matrix.

Note that if $m=1$ then the $E_1$-page is concentrated in columns
$s=0$ and $s=1$ and so collapses at the $E_2$-page for dimensional
reasons. The differential on the $E_1$-page is precisely the boundary
map in the usual long exact sequence. Indeed one can always assemble
such a situation into a long exact sequence.

It is worth commenting that the essential ingredient for the
construction of the spectral sequence is the {\em cube} construction
of link homology, not the particular variant of link homology we
choose to consider. Thus for example one may set up similar spectral
sequences in Khovanov-Rozansky homology.

\section{The rational Khovanov homology of $(3,q)$-torus links.}
\begin{figure}[h]
\centerline{
\psfig{figure= 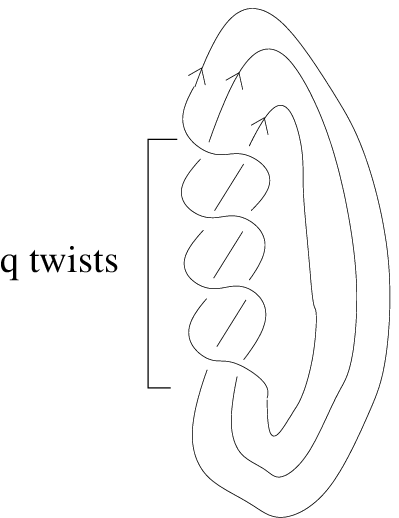}
}
\caption{}
\label{fig:toruslink} 
\end{figure}

In this section we work over $\bQ$ and write $\kh ** D$ for $\khr ** D
{\bQ}$. Our interest is with the torus links $T(3,q)$ which we take to
have negative crossings. We consider the  diagram for $T(3,q)$ obtained as the closure of a three stranded braid as shown in Figure
\ref{fig:toruslink}. When $q$ is a multiple of 3 then $T(3,q)$ is a
three component link, otherwise $T(3,q)$ is a knot.

\begin{thm}\label{thm:torus} Let $N$ be an integer, $N\geq 1$.

(i) The rational Khovanov homology of the $(3,3N)$-torus link is given in Figure \ref{fig:KH33N}.

\begin{figure}[h]
\centerline{
\psfig{figure= 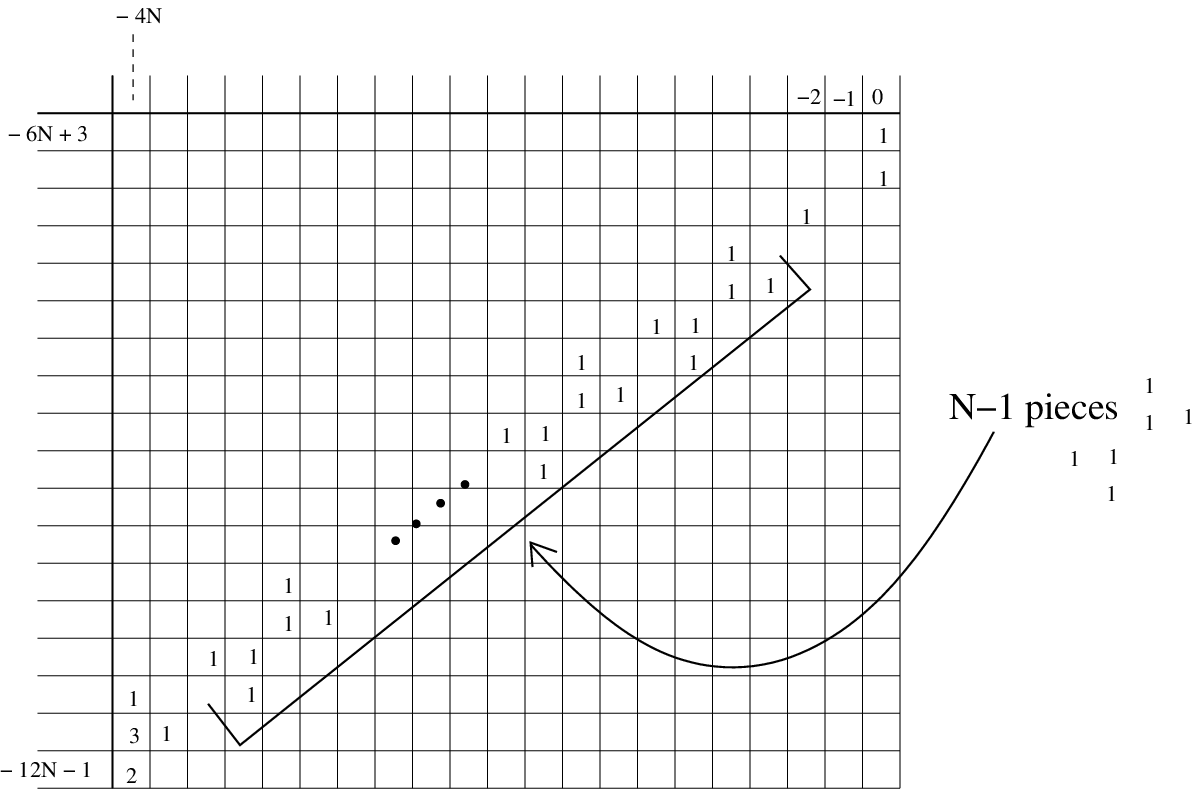}
}
\caption{}
\label{fig:KH33N} 
\end{figure}

\newpage

(ii) The rational Khovanov homology of the $(3,3N+1)$-torus knot is given in Figure \ref{fig:KH33Nplus1}.

\begin{figure}[h]
\centerline{
\psfig{figure= 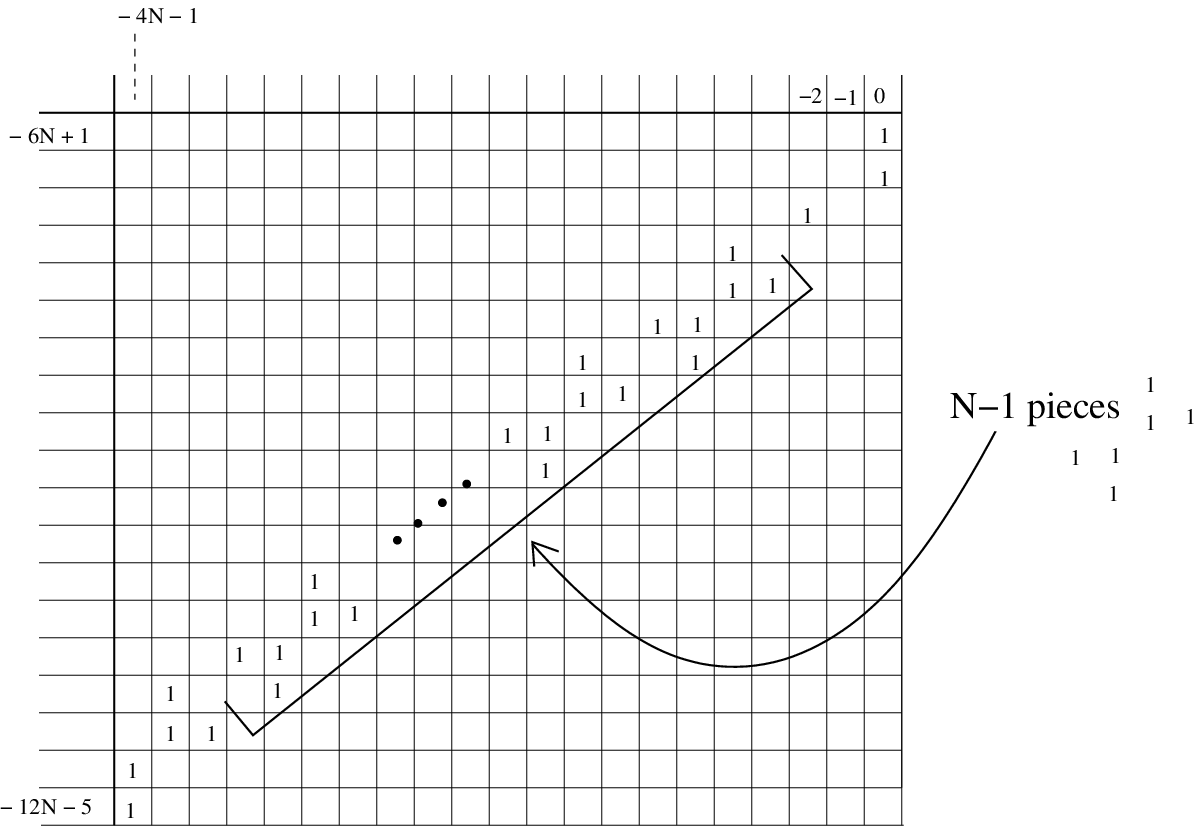}
}
\caption{}
\label{fig:KH33Nplus1} 
\end{figure}

(iii) The rational Khovanov homology of the $(3,3N-1)$-torus knot is given in Figure \ref{fig:KH33Nminus1}.

\begin{figure}[h]
\centerline{
\psfig{figure= 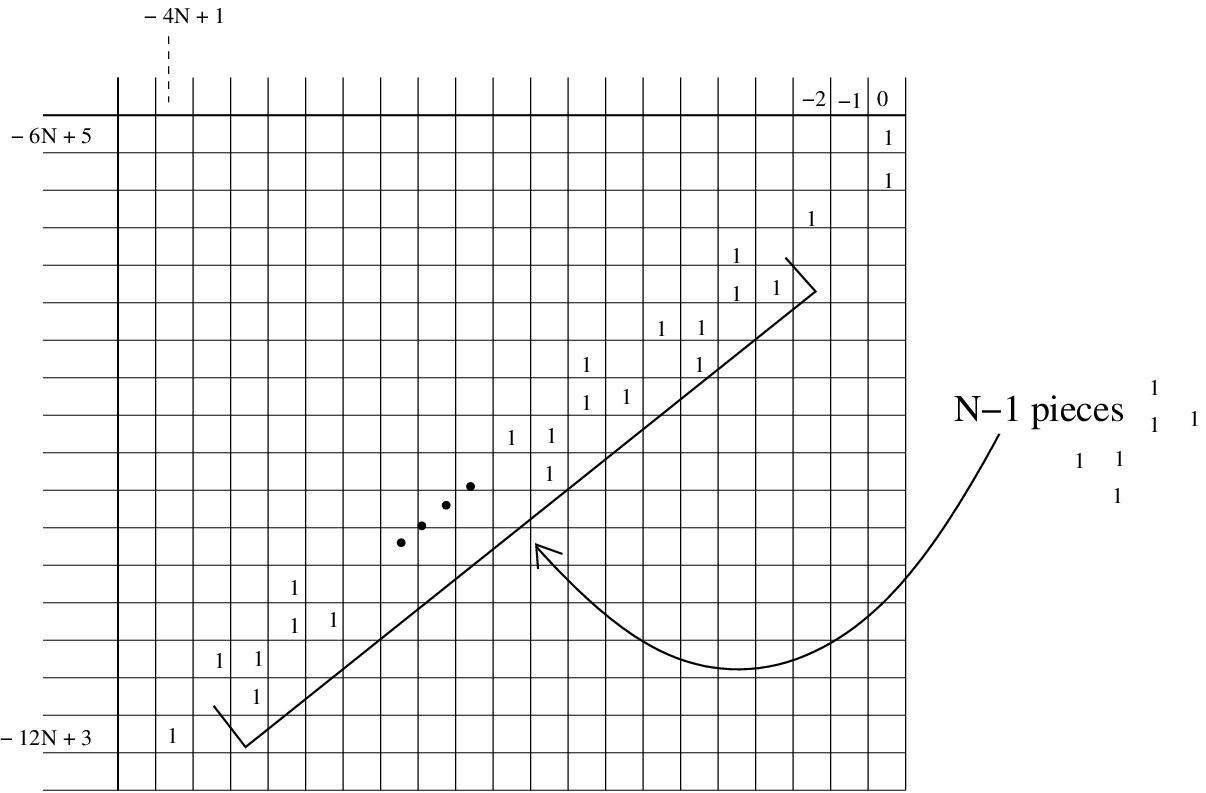}
}
\caption{}
\label{fig:KH33Nminus1} 
\end{figure}

\end{thm}

This has the following corollary.

\begin{cor}
The rational Khovanov homology of the torus links $T(3,3N)$, $T(3,3N+1)$ and $T(3,3N+2)$ occupy exactly $N+2$ diagonals.
\end{cor}

Before proving the theorem we need to make some recollections about
another spectral sequence defined by Lee \cite{lee} (see also
\cite{rasmussen}). Recall that Lee theory 
is a variant of rational Khovanov homology obtained from the same underlying
vector spaces but using a different differential (based on a
different Frobenius algebra). Lee theory is a singly graded theory and
we denote it by $\lee * L$. We summarise the results we need about Lee
theory in the following proposition.

\begin{prop} Let $L$ be an oriented link with $k$ components $L_1, L_2, \ldots , L_k$.

(i) dim$(\lee * L)=2^k$. 

(ii) For every orientation $\theta$ of $L$
there is a generator of homology in degree
\[
 2 \times \sum_{l\in E,m\in \ob{E}} \lk(L_l,L_m)
\]
where $E\subset \{1, 2, \cdots , k \}$ indexes the set of components
of $L$ whose original orientation needs to be reversed to get the
orientation $\theta$ and $\ob{E}=\{1,\ldots,k\}\backslash E$. The
linking numbers $\displaystyle{\lk(L_l,L_m)}$ are the linking number
(for the original orientation) between component $L_l$ and $L_m$.

(iii) There is a spectral sequence converging to $\lee * L$ with $E_1^{s,t} =
\kh {s+t} {4s} L$.
\end{prop}

If we index the $E_1$-page of the spectral sequence by the usual
indexing of the Khovanov homology (rather than of the spectral
sequence) we note that the differential is of
bi-degree $(1,4)$. Indexing the $E_i$-page similarly, the differential
has bi-degree $(1,4i)$.

We note that for a knot Lee theory has two generators in degree
zero. For the $(3,3N)$-torus link (a three component link) Lee theory
has two generators in degree zero and six generators in degree $-4N$.

\begin{proof}{\bf (of Theorem \ref{thm:torus})}
The proof consists of three claims:\\
{ Claim 1:} if the result is true for $T(3,3N-1)$ then the result is true for $T(3,3N)$.\\
{ Claim 2:} if the result is true for $T(3,3N)$ then the result is true for $T(3,3N+1)$.\\
{ Claim 3:} if the result is true for $T(3,3N+1)$ then the result is true for $T(3,3(N+1)-1)$.

Each claim is proved by the same technique, namely we use the
$E_1$-page of the spectral sequence defined in Section 2 to produce
some generators and also to produce some additional {\em possible}
generators. It is not sure that the possible generators are in fact
generators because there may be higher differentials in the spectral
sequence killing them. We then use Lee's spectral sequence to
determine whether or not these possible generators are killed or
not. By playing off one spectral sequence against another in this way,
we do not actually have to explicitly compute any differentials in
either spectral sequence.

{\bf Proof of Claim 1} We will calculate the Khovanov homology of the
link $T(3,3N)$ under the assumption that the Khovanov homology of
$T(3, 3N-1)$ is as given in the statement of the theorem. Consider the
set of crossings consisting of the two top crossings in the braid
diagram (so $m=2$). We have diagrams as presented in Figure
\ref{fig:diags1}.

\begin{figure}[h]
\centerline{
\psfig{figure= 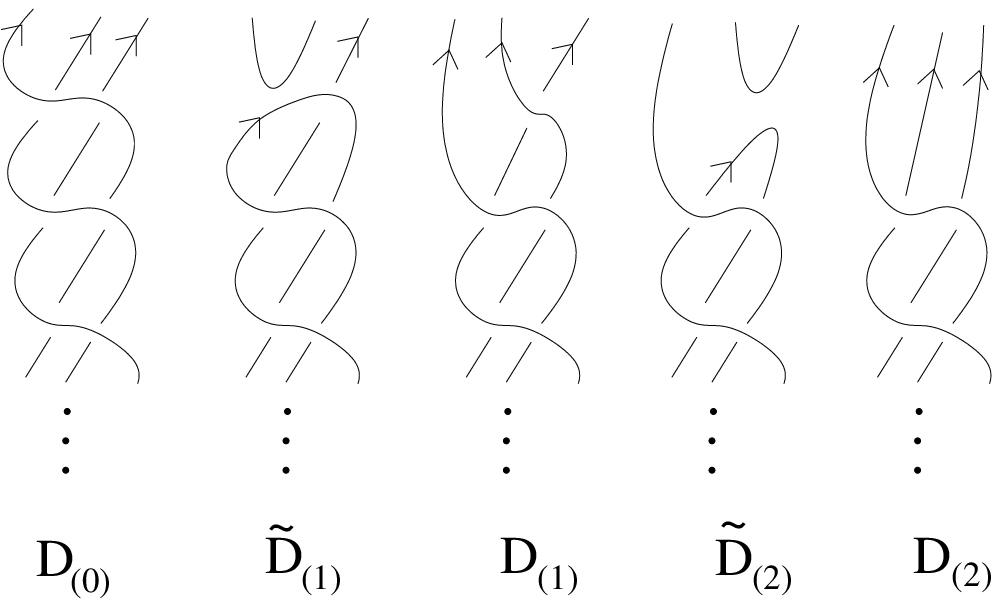}
}
\caption{}
\label{fig:diags1} 
\end{figure}

Note that $\dk 2 = T(3,3N-1)$ and it is easy to see that $\dkb 1 \sim U
\sqcup U$ and $\dkb 2 \sim U$, where $U$ is the unknot. Using the orientations shown in Figure \ref{fig:diags1} one computes
\[
\ol{n}_1^+ = 4n-1, \;\; \ol{n}_1^- = 2N, \;\;
\ol{a}_1 = 4N,\;\; \ol{b}_1=12N-1
\]
and
\[
\ol{n}_2^+ = 4n-1, \;\; \ol{n}_2^- = 2N-1, \;\;
\ol{a}_2 = 4N,\;\; \ol{b}_2=12N-1.
\]
From Proposition \ref{prop:ss} we have
\[
E_1^{0,t} = \kh {t+4N}{j+12N-1}{U\sqcup U},
\]
\[
E_1^{1,t} = \kh {t+4N+1}{j+12N}{U},
\]
\[
E_1^{2,t} = \kh {t+2}{j+2}{T(3,3N-1)}.
\]

When $s=0$ we see that  $E_1^{0,t}=0$ unless $t=-4N$ and $j=-12N-1$, $j=-12N +1$ or $j=-12N+3$. Similarly, $E_1^{1,t}=0$ unless $t=-4N-1$ and  $j=-12N-1$ or 
$j=-12N +1$ and $E_1^{2,t}$ is zero unless $-12N+1 \leq j \leq -6N +3$.

For $j>-12N+3$ the $E_1$-term of the spectral sequence is concentrated
in the column $s=2$ and hence collapses for dimensional
reasons. Thus,
\[
\kh i j {T(3,3N)} \cong E_1^{2,i-2} = \kh i {j+2} {T(3,3N-1)}.
\]
We need to consider the three cases $j=-12N-1$, $j=-12N+1$ and
$j=-12N+3$. The $E_1$ pages are give in Figure \ref{fig:E1claim1}.

\begin{figure}[h]
\centerline{
\psfig{figure= 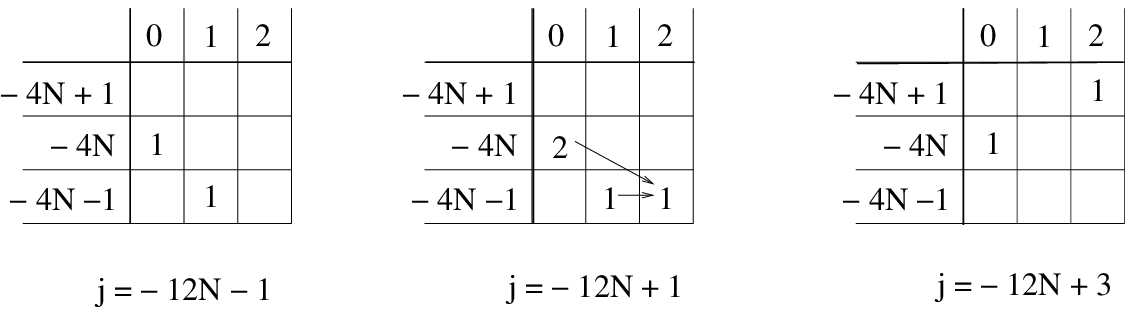}
}
\caption{}
\label{fig:E1claim1} 
\end{figure}

For $j=-12N-1$ and $j=-12N+3$ there are no differentials for
dimensional reasons thus the spectral sequence collapses at $E_1$. For
$j=-12N+1$ there is a possible $d_1$ and a possible $d_2$ (but not
both) as shown in Figure
\ref{fig:E1claim1}. Thus for $T(3,3N)$ we have the situation presented in 
Figure \ref{fig:claim1}, where {\em possible} generators are circled.

\begin{figure}[h]
\centerline{
\psfig{figure= 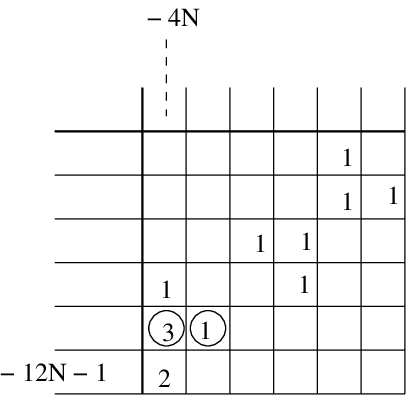}
}
\caption{}
\label{fig:claim1} 
\end{figure}

The three possible generators in bi-degree $(-4N, -12N+1)$ must all
indeed be generators because we require at least six generators in
homological degree $-4N$. This is because Lee theory in this degree
has six generators and due to Lee's spectral sequence these must show up
in Khovanov homology. 

The possible generator in bi-degree $(-4N+1,-12N+1)$ is also a
generator. If we look at the $E_1$ page for $j=-12N+1$ then since the
three generators on the line $s+t=-4N$ survive until $E_\infty$ (as
shown in the previous paragraph) then there is nothing to kill the
remaining generator. (Alternatively, the generator in bi-degree
$(-4N+2,-12N+5)$ must be killed in Lee's spectral sequence and the
only possible way this can happen is for the possible generator in
bi-degree $(-4N+1,-12N+1)$ to be present. To see this, recall that
indexed this way the differential $d_i$ in Lee's spectral sequence has
bi-degree $(1,4i)$.)

Thus we end up computing $\kh * * {T(3,3N)}$ as presented in the theorem.

{\bf Proof of Claim 2} Consider the link $T(3,3N+1)$ and as above take
the set of crossings to be the two top crossings in the braid
diagram. We have diagrams as presented in Figure \ref{fig:diags2}.

\begin{figure}[h]
\centerline{
\psfig{figure= 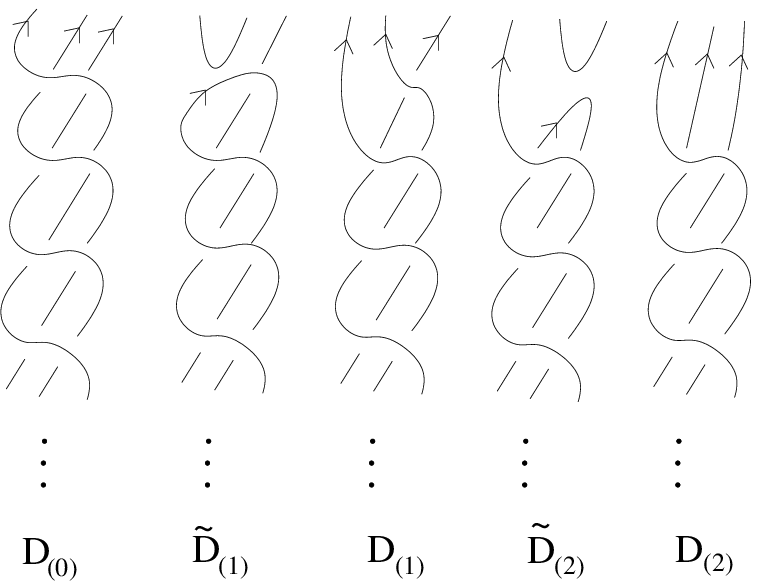}
}
\caption{}
\label{fig:diags2} 
\end{figure}

Note that $\dk 2 = T(3,3N)$ and it is easy to show $\dkb 1 \sim U
$ and $\dkb 2 \sim U \sqcup U$. Using the orientations shown in Figure \ref{fig:diags2} one computes
\[
\ol{n}_1^+ = 4N, \;\; \ol{n}_1^- = 2N+1, \;\;
\ol{a}_1 = 4N+1,\;\; \ol{b}_1=12N+2
\]
and
\[
\ol{n}_2^+ = 4N, \;\; \ol{n}_2^- = 2N, \;\;
\ol{a}_2 = 4N+1,\;\; \ol{b}_2=12N+2.
\]
Thus we have
\[
E_1^{0,t} = \kh {t+4N+1}{j+12N+2}{U}
\]
\[
E_1^{1,t} = \kh {t+4N+2}{j+12N+3}{U\sqcup U}
\]
\[
E_1^{2,t} = \kh {t+2}{j+2}{T(3,3N)}
\]

For $s=0$ we must have $j$ in the range $-12N-3\leq j \leq -12N-1$,
for $s=1$ in the range $-12N-5\leq j \leq -12N-1$ and for $s=2$ in the
range $-12N-3\leq j \leq -6N+1$. For $j> -12N-1$, as in the previous
case, we instantly see that the result is as claimed. For the three
remaining $j$-values we have $E_1$-pages as given in Figure
\ref{fig:E1claim2}.

\begin{figure}[h]
\centerline{
\psfig{figure= 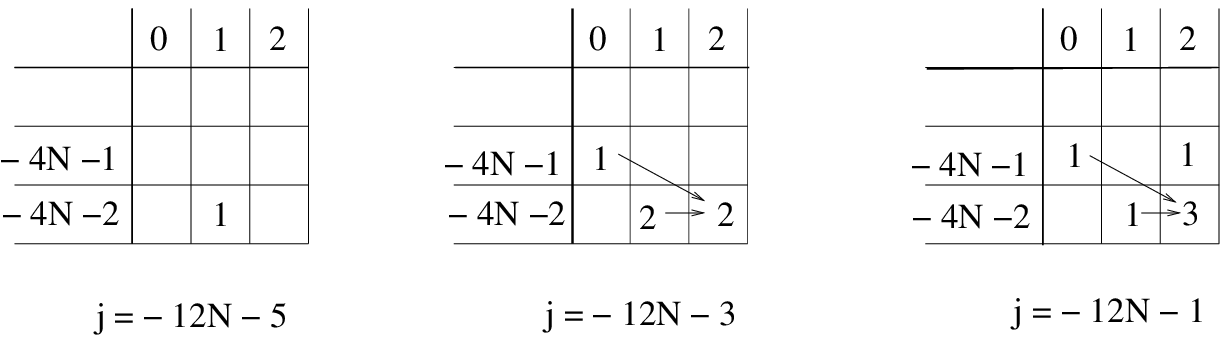}
}
\caption{}
\label{fig:E1claim2} 
\end{figure}

For $j=-12N-5$ there are no differentials for dimensional reasons, but for
$j=-12N-3$ and $j=-12N-1$ there are possible differentials.  The
situation is presented in
Figure \ref{fig:claim2}, where, as above, {\em possible} generators are circled.

\begin{figure}[h]
\centerline{
\psfig{figure= 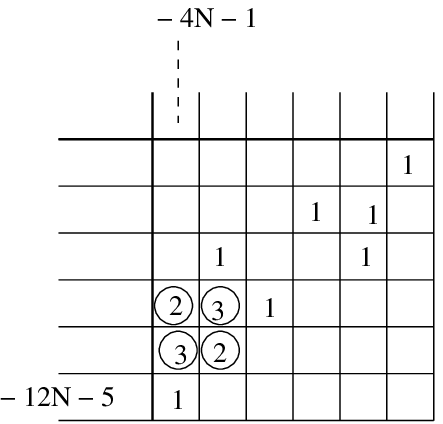}
}
\caption{}
\label{fig:claim2} 
\end{figure}

Consider the two possible generators in bi-degree $(-4N, -12N-3)$. There
cannot be generators in this bi-degree since they would
appear in the $E_\infty$-page of Lee's spectral sequence. However, 
$T(3,3N+1)$ is a knot so the $E_\infty$-page has only two generators and
these lie on the line $s+t=0$.

Now look at the $E_1$-page for $j=-12N-3$. We have just argued that
the two generators on the line $s+t=-4N$ must be
killed. There are two possible ways this might happen, but either way
one is left with one generator on the line $s+t=-4N -1$ and this
must survive to $E_\infty$.

A similar argument holds for the two possible generators in bi-degree
$(-4N-1, -12N-1)$ and one is left with one generator in
bi-degree $(-4N, -12N-1)$.

{\bf Proof of Claim 3} This is very similar to the previous arguments
so we present this case only briefly. We follow the same orientation 
convention as above for the diagrams. We have
\[
E_1^{0,t} = \kh {t+4N+3}{j+12N+8}{U},
\]
\[
E_1^{1,t} = \kh {t+4N+3}{j+12N+6}{U},
\]
\[
E_1^{2,t} = \kh {t+2}{j+2}{T(3,3N+1)}.
\]

For $j> -12N-5$ there is nothing to do and for the remaining
$j$-values of interest we have $E_1$-pages as given in Figure
\ref{fig:E1claim3} leading to the generators and possible generators presented in Figure \ref{fig:claim3}.

\begin{figure}[h]
\centerline{
\psfig{figure= 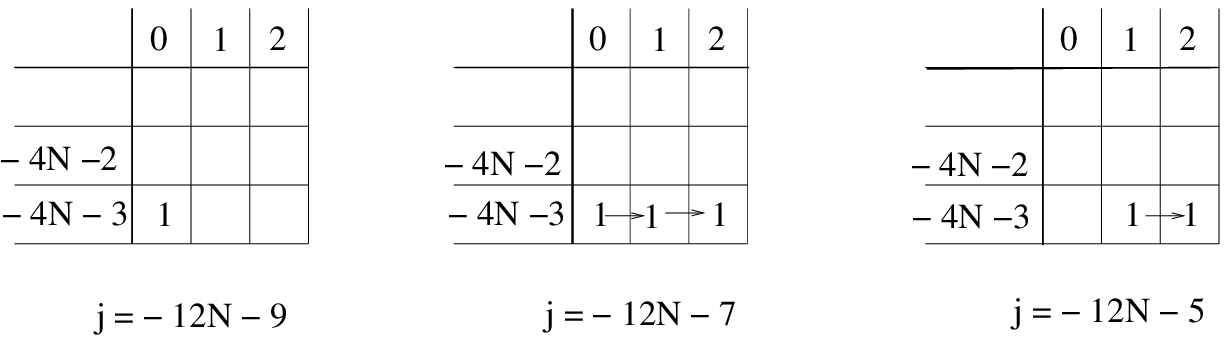}
}
\caption{}
\label{fig:E1claim3} 
\end{figure}

\begin{figure}[h]
\centerline{
\psfig{figure= 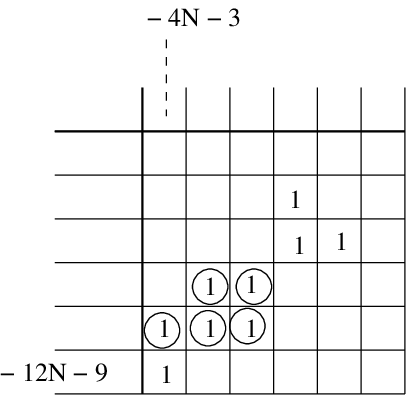}
}
\caption{}
\label{fig:claim3} 
\end{figure}

As above it is easy to see that the possible generator in bi-degree
$(-4N-3, -12N-7)$ cannot survive to $E_\infty$ in Lee's spectral
sequence so must be killed. There is only one
possibility which leaves one generator in homological
degree $-4N-1$. When $j=-12N-5$ the two generators both survive
because they are needed in Lee's spectral sequence to kill the
generators in bi-degree $(-4N-3,-12N-9)$ and $(-4N,-12N-1)$.

Finally we note that the inductive process of the above three claims
has a beginning because the Khovanov homology of $T(3,2)$ is easily
calculated (even by hand), and $T(3,3)$, $T(3,4)$ and $T(3,5)$ can also
be computed (by computer or using the spectral sequence - the
computations are similar, though not identical, to those above). These
cases are seen to have the required form.
\end{proof}

The Khovanov homology of positive crossing $(3,q)$-torus links can be
computed from the above by recalling that the rational Khovanov
homology of the mirror image $L^!$ of a link $L$ can be computed as
$\kh i j {L^!} = \kh {-i}{-j} L$.

\section*{Acknowledgements}
Thanks to M. Mackaay and J. Rasmussen for comments on a draft version.

\end{document}